\documentclass[a4paper,11pt,reqno]{amsart}
\usepackage{amsmath,amssymb,amsthm,enumerate}

\usepackage[latin1]{inputenc}
\usepackage[english]{babel}
\newcommand{\hilb}{{\mathcal{H}\textnormal{ilb}}}
\newcommand{\hilbp}{{\mathcal{H}\textnormal{ilb}_{p(z)}^n}}

\usepackage{amscd}
\newtheorem{Theorem}{\quad Theorem}[section] 

\newtheorem{Definition}[Theorem]{\quad Definition} 

\newtheorem{Proposition}[Theorem]{\quad Proposition} 

\newtheorem{Corollary}[Theorem]{\quad Corollary} 

\newtheorem{Lemma}[Theorem]{\quad Lemma} 

\newtheorem{Example}[Theorem]{\quad Example}

\newtheorem{remark}[Theorem]{\quad Remark} 

\newtheorem{procedure}[Theorem]{\quad Procedure} 

\newcommand{\In}{\textnormal{In}} 
\newcommand{\LM}{\textnormal{LM}} 
 \newcommand{\ed}{\textnormal{ed}}   
 \newcommand{\lcm}{\textnormal{lcm}} 
 \newcommand{\tail}{\textnormal{tail}} 
  \newcommand{\htail}{\textnormal{htail}} 

\newcommand{\Proj}{\textnormal{Proj}\,}

\newcommand{\TT}{\mathbb{T}}
\newcommand{\ZZ}{\mathbb{Z}}
\newcommand{\Af}{\mathbb{A}}
\newcommand{\PP}{\mathbb{P}}
\newcommand{\Tx}{\mathbb{T}_{X}}
\newcommand{\Txc}{\mathbb{T}_{X,C}}
\newcommand{\Tc}{\mathbb{T}_{C}}
\newcommand{\Ty}{\mathbb{T}_{Y}}
\newcommand{\Txy}{\mathbb{T}_{X,Y}}

\newcommand{\Gr}{Gr\"obner\, }
\newcommand{\idR}{\mathcal{R}}
\newcommand{\n}{\mathbf{n}}

\newcommand{\m}{\mathbf{m}}
\newcommand{\vv}{\mathbf{v}}

\newcommand{\p}{\mathbf{p}}
\newcommand{\q}{\mathbf{q}}

\newcommand{\Aid}{\mathcal{A}}
\newcommand{\St}{\mathcal{S}t}
\newcommand{\htt}{\mathrm{ht}}

\date{}

\title{IDEALS WITH AN  ASSIGNED INITIAL IDEAL}

\author{Margherita Roggero and Lea Terracini} 

\begin{document}

\footnotetext{Written with the support of the University Ministry funds.} 
\begin{abstract} 
The stratum $\St(J,\prec)$ (the homogeneous stratum $\St_h(J,\prec) $ respectively) of a monomial ideal $J$ in a polynomial ring $R$  is the family of all  (homogeneous) ideals   of   $R$  whose initial ideal with respect to the   term order $\prec$ is $J$.   $\St(J,\prec)$ and  $\St_h(J,\prec)$ have a natural structure of affine schemes. Moreover they are  homogeneous   w.r.t.  a non-standard graduation called level. This property allows us to draw  consequences that are interesting from both a theoretical and a computational point of view. For instance  a smooth stratum is always isomorphic to an affine space (Corollary \ref{usolivelli}). As applications, in \S \ref{sec:lex}   we prove that   strata and homogeneous strata w.r.t. any term ordering $\prec$ of every  saturated \texttt{Lex}-segment ideal $J$ are smooth. For $\St_h(J,\texttt{Lex})$ we also give a formula for the dimension. In the same way in \S \ref{sec:revlex} we take in consideration   any ideal $\idR$   in $k[x_0, \dots, x_n]$  generated by a saturated \texttt{RevLex}-segment ideal in $k[x_0,x_1,x_2]$. We also prove that   $\St_h(\idR,\texttt{RevLex})$ is smooth and give a formula for its  dimension.
 \end{abstract}

\maketitle
\section{Introduction} Let us consider a polynomial ring $R$  over a field $k$, a term order $\prec$ on $R$ and a monomial ideal $J\subset R$. In this paper we study the family $\St(J,\prec)$ of all ideals in $R$ having $J$ as initial ideal with respect to $\prec$ in order to provide it  with a structure of algebraic scheme, when possible. Sometimes, this family can be too huge to allowe such a structure (see Example \ref{grosso}). 
In order to avoid this situation  we  either  consider the family $\St_h(J,\prec))$ of homogeneous ideals or  pose some extra condition on the term order:   we say that $\prec$ is \emph{{reliable}} if  for every monomial $\mathbf{m}$  the set of monomials  $\mathbf{n}\prec \mathbf{m}$  is finite. Under   one of the previous conditions, the family of ideals $I$ such that $\In(I)=J$ is  in a natural way an affine scheme, that we call  \emph{{stratum}} (\emph{{homogeneous stratum}} respectively) following the terminology introduced in \cite{NS}, where homogeneous strata with respect to \texttt{RevLex} are related with a suitable stratification of Hilbert schemes.

 There is a very natural way for explicitly building on the algebraic structure of  strata,  making use of Buchberger's algorithm  for \Gr bases.   In this way the ideal $J$ corresponds to the origin in both $\St(J,\prec)$ and  $\St_h(J,\prec)$ and moreover   the set-theoretic inclusion $\St_h(J,\prec)\subset \St(J,\prec)$  is an algebraic map, because $\St_h(J,\prec)$ becomes   the  section of   $\St(J,\prec)$  with  a suitable  linear space.  

A possible   objection  to our construction   is that of giving rise to geometrical objects having a well defined support, but possibly   different structures of affine schemes, because   the procedure of reduction with respect to  polynomials that are not a \Gr basis is not unique.    This  is claimed for instance by Robbiano in   \cite[\S 3]{Ro}.   Even though this remark can  look reasonable at first, however we want to emphasize that it is  not correct and that strata  have a unique and well defined structure of affine schemes that does not depend on the procedure used   to obtain it and that can be reducible or not reduced. For a proof of the well definiteness of strata we refer to     \cite{LR}. Reducible or non-reduced strata are presented in Example \ref{riducibile} and in \S \ref{sec:esempi}.

There are no difficulties to give an implementation of our construction  in order to obtain explicitly the ideal of a stratum in a suitable polynomial ring, that is to realize the stratum as a closed subscheme of an affine space $\Af^N(k)$. However, this is computationally very heavy, because it requires to introduce a very big number of variables and the ideal of the stratum is in general given by a huge number  of generators, so that even a basic  study of strata, like for instance  dimensions or singularities, requires   times of calculation absolutely unsatisfactory or even not affordable from a practical point of view.

In this paper we present a refined method for the construction of strata, that allows  to  reduce in a sensible way the time of calculation and the amount of memory necessary to compute   equations of a stratum. Moreover  our method leads to   some interesting theoretical consequences.

The main idea is the following:
the stratum is given by an   ideal $\Aid(J)$ in a suitable polynomial ring, which in general is not  homogeneous  w.r.t. the usual graduation; however we can define a new graduation, the \emph{level}, such that  $\Aid(J)$ turns out to be level-homogeneous.

Indeed, the levels allow us to prove the following result (Corollary \ref{usolivelli}):

\medskip

\centerline{\textit{$J$ is a smooth point of $\St(J,\prec)$    $\Longleftrightarrow$  $\St(J,\prec)\cong \Af^s(k)$}.}

\medskip

The same holds for $\St_h(J,\prec)$. This fact is claimed (but not proved) in \cite{NS}. The fact that strata are level-homogeneous has a more general consequence: every stratum can be isomorphically embedded in its  Zariski tangent space at the origin, whose dimension is in fact the embedding dimension of the stratum. The level-homogeneity of strata is the key point of   Procedure \ref{procedura} that gives an algorithm for a simplified construction of a set of equations for the strata. A Maple12 implementation of Procedure \ref{procedura}  is available at: 

\noindent{\texttt{http://www2.dm.unito.it/paginepersonali/roggero/InitialIdeal(Maple)/}}

In the last two  sections, to the aim of illustrate the potential of levels, we present the direct computation of strata in two important families of ideals.

In \S \ref{sec:lex}   we consider the case of the saturated \texttt{Lex}-segment ideals $J$ that are of the utmost importance in the theory of Hilbert schemes. As well known, every Hilbert scheme $\hilbp$ parameterizing subschemes in $\PP^n$ contains one and only one point corresponding to such an ideal (see \cite{M}). This point, the \texttt{Lex}-point,  play a key role in most of the main general results on Hilbert schemes. For instance, in \cite{RS}  it is proved that  $J$ belongs to only one irreducible component of $\hilbp$ (called  \texttt{Lex}-component or the Reeves and Stillman component) and that it is a smooth point on it. By the universal properties of Hilbert schemes,    the homogeneous stratum of $J$ can be embedded as a locally closed subscheme of the \texttt{Lex}-component.
In \cite{NS} it is claimed  that the smoothness of   $\St_h(J,\texttt{RevLex})$ at $J$ can be deduced   from that of the \texttt{Lex}-component. This argument clearly holds  when  $\St_h(J,\texttt{RevLex})$ corresponds to an open subset of $\hilbp$ and this is in general not true (see Remark \ref{ultima}).

However, we can give a direct proof of the smoothness of  both strata and homogeneous strata w.r.t. any term ordering $\prec$ of every  saturated \texttt{Lex}-segment ideal $J$. In particular   we   give a formula for the dimension of $\St_h(J,\texttt{Lex})$ and  a   condition for $\St_h(J,\texttt{Lex})$  to be an open subset of $\hilbp$. 

Using the same technique, in \S \ref{sec:revlex} we   consider     any ideal $\idR$   in $k[x_0, \dots, x_n]$  generated by a saturated \texttt{RevLex}-segment ideal in $k[x_0,x_1,x_2]$. If $n=2$, these ideals $\idR$ are the generic initial ideals w.r.t. \texttt{RevLex} of the ideals of sets of general points in $\PP^2$. For  $n\geq 3$ the ideals $\idR$ defines  arithmetically Cohen-Macaulay subschemes of codimension 2 in $\PP^n$. Again, we are able to  prove that  the homogeneous strata $\St_h(\idR,\texttt{RevLex})$ are smooth, hence  isomorphic to affine spaces,   and to give a formula for their dimension.
  
 \section{General settings and definition  of strata}\label{Sec:stratum} 

We shall denote   by $ X $ the set of variables  $x_1, \dots, x_n$ and by $k[ X ]$ the polynomial ring $k[x_1, \dots, x_n ]$ over a field $k$. In the same way,   $C $ will be the set of variables $ c_{i,\alpha } $ (that  will be introduced in the next section) and   $X,C $ the union of the two sets of variables.
 $\Tx$ will be the semigroup of monomials  in the variables  $X$; analogous meanings will have  $\Tc$ and $\Txc$.   A monomial $x_1^{\alpha_1}\cdots x_n^{\alpha_n}\in \Tx$ will be written as $X^\alpha $ where $\alpha=[\alpha_1, \dots,\alpha_n]$ is the ordered list of exponents. If $\prec$ is any term ordering on $\Tx$, we will also denote by $\prec$ the induced total ordering on $\ZZ^n$; in this way $(\ZZ^n,+, \prec)$ is an ordered group. 
 \begin{Definition} We shall call \emph{tail} and  \emph{homogeneous tail} of  $X^\alpha \in \Tx$ w.r.t.  $\prec$, the sets $\tail(X^\alpha)=\{ X^\beta \in \Tx / \beta \prec \alpha \}$ and  $\htail(X^\alpha)=\{ X^\beta \in \Tx / \beta \prec \alpha \hbox{ and } \vert \beta\vert=\vert \alpha \vert \}$ respectively. If $\tau$ is any subset of $k[X]$, $\tail_\tau (X^\alpha)$ and $\htail_\tau(X^\alpha)$ will be $\tail(X^\alpha)\setminus \tau$ and $\htail(X^\alpha)\setminus \tau$ respectively .
 \end{Definition}
 Of course, $\htail(X^\alpha)$ is a finite set, while $\tail(X^\alpha)$ may have infinitely many elements. 
\begin{Example}\label{grosso} Consider $k[x,y]$ equipped by the term order \texttt{Lex} with $y\prec x$ and the monomial ideal $J=(x)$. The family of ideals $I$ having $J$ as initial ideal contains   all ideals of the form  $(x+g(y))$ with $g$ varying in $k[y]$ and therefore it depends on infinitely many free parameters.
\end{Example}
A \emph{reliable term order} $\prec$ in $\Tx$ is a term order such that  every monomial    has a finite tail. For instance, every graded term order is reliable. More generally, if $\prec$ is defined using a matrix with integer entries $M$, we can check if it is reliable through the criterion given by following :
\begin{Lemma}$\prec$ is reliable $\Longleftrightarrow$ the first row of $M$ has strictly positive entries.
\end{Lemma}
\begin{proof} If the first row $[a_1, \dots , a_n]$ has strictly positive entries,  $X^\beta \prec X^\alpha$ only if $\sum \beta_i a_i \leq \sum \alpha_i a_i$, where $\alpha_i a_i\geq 0$, hence each $\beta_i$ is upper bounded by $a/a_i$. On the converse, if we have $a_i=0$ for some $i$, and of course also $a_j>0$ for some $j$, then   $\tail(x_j)$ contains all the monomials $x_i^r$ for every $r> 0$.
\end{proof}
When general tails are concerned, we will always assume that the term ordering is reliable.  

If $F$ is a polynomial, $\LM_\prec (F)$, or for short $\LM(F)$,  is the leading monomial of $F$ with respect to $\prec$. If $I$ is an ideal, $\In_\prec (I)$, or for short $\In(I)$, will be the initial ideal of $I$ w.r.t. $\prec$ and $B_I$ will be its reduced \Gr basis. If $I$ is monomial, then of course $\In(I)=I$ and $B_I$ is its  monomial basis.

 In this paper we take into consideration some families of ideals $I$ in $k[X]$ having the same initial ideal $J$.   If $B_J=\{ X^{\alpha_1}, \dots, X^{\alpha_s}\}$,  the reduced \Gr basis $B_I$ of every such ideal $I$  is formed   by $s$ polynomials  of the type $X^{\alpha_i} + \sum a_{i \beta}X^\beta $, where $a_{i\beta}\in k$ and  $X^\beta \in \tail_J (X^{\alpha_i})$.  
Let us consider the polynomials  
\begin{equation}\label{polF}  F_i=X^{\alpha_i} + \sum c_{i \alpha}X^\beta
\end{equation} of the previous type, but whose  coefficients $c_{i \beta}$ are new variables and  let  $\mathcal{I}$ be the ideal in $k[X, C ]$ generated by the set $\mathcal{B}=\{F_1, \dots, F_s\}$.
 The reduced \Gr basis of $I$   can be obtained from $\mathcal{B}$ by specializing in $k$ the  parameters $C$.  However, not every specialization of $ \mathcal{B}$ returns a \Gr basis and so an ideal $I$ with $\In(I)=J$.  In order to characterize  the good specializations, we consider the following construction. 
 
  Let us fix any term order $\prec_1$ on $\TT_C$ and let $\prec'$ the term order on $\TT_{X,C}$ which is an elimination order for the variables $X$ and  coincides with $\prec$ and $\prec_1$  on $\TT_{X}$ and $\TT_{C}$ respectively. 
  Moreover let us consider  some reductions $H_{ij}$ by ${\mathcal{B}}$ of the  $S$-polynomials of each pair of elements in ${\mathcal{B}}$: notice that this procedure commutes with the specialization in $k$ of the variables $C$, because the leading monomials  of  the $F_i$'s  are monomials in the variables  $X$. A specialization of ${\mathcal{B}}$ is a \Gr basis if and only if all the $H_{ij}$ specialize to  zero.  This construction not only gives a characterization of all the ideals $I$ such that $\In(I)=J$, but also leads us to obtain this family as the set of closed points of an affine subscheme of  the affine space $\Af^N$ over the field $k$,  where $N=\sum_i   \vert \tail_J(X^{\alpha_i}) \vert$.

\begin{Definition}\label{def:strata} Let $J\subset k[X]$ be a monomial ideal, $B_J=\{ X^{\alpha_1}, \dots,X^{\alpha_s}\}$ its basis. We will call  \emph{stratum of} $J$ w.r.t. $\prec$ the    affine subscheme $\St(J,\prec)$ of $\Af^N$ defined by the ideal $\mathcal{A}(J,\prec)$ in $k[C]$ generated by all the coefficients  in $k[C]$ of monomials in the variables $X$ appearing in a set of  polynomials $\{ H_{ij}, \ 1\leq i < j \leq s\}$. The \emph{homogeneous stratum} $\St_h(J, \prec)$ is the affine scheme obtained in the same way but using   homogeneous tails instead of tails.

More generally,  if $\tau=[\tau_1, \dots, \tau_s]$ is any list of $s$ sets of monomials, we will denote by $\St (J  /  \tau,\prec)$ and by $\mathcal{A}(J/\tau,\prec)$ the subscheme and its ideal obtained as before substituting,  for every $i=1, \dots, s$,  $\tail_{\tau_i\cup J}(X^{\alpha_i})$ to  $\tail_J(X^{\alpha_i})$. When the term order $\prec$ is clearly fixed, we will omit it in the symbols   and write simply $\St(J)$,  $\mathcal{A}(J)$ and so on.
\end{Definition}
It is easy to see that $\St (J / \tau )$ is   the scheme-theoretical intersection of $\St (J )$  with the linear space defined by $c_{i\beta}=0$, where  for every $ i=1, \dots, s$  $ X^\beta $ belongs to $ \tail_{ J}(X^{\alpha_i})\cap \tau_i$. Of course all these objects are non-empty, because they contain  at least a closed point, the origin, corresponding to the trivial case $I=J$. 

\begin{Example} \label{escluso} Let $J$ and $J'$ be monomial ideals in $k[X]$ and $\prec$ any term ordering. We can consider the  subscheme   of  $\St (J)$ whose   closed points    correspond to    ideals  $I\in \St(J)$ such that their reduced \Gr bases are    also reduced modulo $J'$. Of course this subscheme is empty if some monomial in $J'$ divides a monomial in the basis of $J$.  If this does not happen,  this subscheme is  $\St(J /  \tau)$, where $\tau_i= J'$ for every $i$. We will denote it by $\St(J /  J')$.
\end{Example}

\begin{Example} \label{es:omo} Let $B_J=\{ X^{\alpha_1}, \dots,X^{\alpha_s}\}$ be the basis of the monomial ideal $J$ in $k[X]$. If $\tau$ is the list of subsets $\tau_i=\oplus_{i\neq \vert \alpha_i \vert} k[X]_i$, then $\St(J / \tau)=\St_h(J)$.
\end{Example}

We do not present at this point many explicit examples of strata, because in most non-trivial cases the ideal  $\mathcal{A}(J/\tau)$ lives in the polynomial ring  $k[C]$ which has  in general \lq\lq a very big\rq\rq\ number of variables. In the next section we  will see how to obtain a new ideal  $\tilde{\Aid}(J/\tau)$ for $\St(J /  \tau)$ in a polynomial ring    $k[C']$, where $C'$ is the smallest possible subset of $C$ such that $\St(J /  \tau)$ can be embedded in an affine space of dimension $\vert C' \vert$. However we prefer to present here at least   a couple of (necessarily trivial) examples in order to show that strata can in fact be   reducible and  non-reduced schemes. 

\begin{Example}\label{riducibile} Let $J=(xy,y)\subset R=k[x,y]$ and let $\tau=[R \setminus \{x\}, R]$ and $\tau'=[R\setminus\{ x\}, R\setminus \{y\}]$. In the first case $\mathcal{B}=(xy+cx,y^2)$, the only $S$-polynomial is $y(xy+cx)-x(y^2)=cxy$ that reduces to $H_{12}=-c^2x$. Hence $\St(J /  \tau)$ is a non-reduced subscheme of $\mathbb{A}^1$ because  it is defined by the ideal $(c^2)\subset k[c]$.

 In the second case $\mathcal{B}'=(xy+c_1x,y^2+c_2y)$, the only $S$-polynomial is $y(xy+c_1x)-x(y^2+c_2y)=(c_1-c_2)xy$ that reduces to $H_{12}=-c_1(c_1-c_2)x$. Hence $\St(J /  \tau')$ is a reducible subscheme of $\mathbb{A}^2$ because it is defined by the ideal $(c_1(c_1-c_2))\subset k[c_1,c_2]$.
\end{Example}

\section{Strata are homogeneous varieties}\label{sec4}

In this section we fix  a term ordering $\prec$ on $k[X]$ and a monomial ideal $J\subset k[X]$ with $n=\vert X\vert$ and $B_J=\{X^{\alpha_1}, \dots, X^{\alpha_1}\}$. What we are going to prove in this section holds true for every stratum  $\St(J,\prec)$, $\St_h(J,\prec)$ or more generally  $\St(J /  \tau,\prec )$ and any $\tau$. Then we will simply denote by $\St(J)$  either one of them and by $\mathcal{A}(J)$ its ideal.

\begin{Definition} We shall denote by  $\Lambda$ the grading  of either $k[X, C]$  or $k[C]$ over the totally ordered group $(\ZZ^n,+,\prec)$ given by   $\Lambda(X^{\alpha})=\alpha$ and   $\Lambda(c_{i\beta})=\alpha_i -\beta $. We shall call $\Lambda(c_{i\beta})$ the \emph{level} of $c_{i\beta}$. 
\end{Definition}
 
 As we shall use also the usual grading where all   variables have degree 1, we will   always explicit   the symbol $\Lambda$ when the above defined grading is concerned (so, $\Lambda$-degree $\lambda$ with $\lambda\in \ZZ^n$, $\Lambda$-homogeneous of degree $\lambda$,  $\Lambda$- component of degree $\lambda$ etc.) and leave the simple term when the usual grading is involved (so, degree $r$ with $r\in \ZZ$, homogeneous of degree $r$, component of degree $r$, etc.).
\begin{Lemma}\label{stessolivello} In the above settings:
\begin{itemize}
\item[i)] $\Lambda$ is a positive grading, that is for every $X^\alpha C^\gamma\in \Txc$, $\Lambda(X^\alpha C^\gamma)\succ \Lambda(1)$.
\item[ii)] the ideal  $\Aid (J)$ of the stratum $\St(J)$  is   $\Lambda$-homogeneous.  
\end{itemize}
\end{Lemma} 
 \begin{proof} The first item immediately follows from the definitions. In fact $\Lambda(x_j)\succ \Lambda(1)$ because   $\prec$ is a term ordering and $\Lambda(c_{i\beta})\succ \Lambda(1)$  because $X^\beta \in \tail(X^{\alpha_i})$ and so $X^{\alpha_i}\succ X^\beta$.
 
For the second item, it is sufficient to observe that the  function $\Lambda$ is, by definition, constant over all the monomials of $\Txc$ used in   the polynomials $F_i$ that we use in the construction of $\Aid(J)$.  Hence  the  polynomials $F_i$, their $S$-polynomials $H_{ij}$ and their reductions  are $\Lambda$-homogeneous. Finally if we collect a $\Lambda$-homogeneous polynomial  in $k[X,C]$ with respect to the monomials of  $\Tx$, the coefficients are $\Lambda$-homogeneous polynomials in $k[ C]$.   
 \end{proof}

We now present some properties of $\Lambda$-homogeneous ideals of a polynomial ring graded over an ordered group $(\ZZ^n,+,\prec)$. As we will apply these properties to any ideal $\mathcal{A}\subseteq k[C]$, we denote by $C$ the set of variables, but for  sake of simplicity $C=\{c_1, \dots, c_s\}$ (instead of $C=\{c_{i\beta}, i=1, \ \dots, n, \beta \in \tail(X^{\alpha_i})\}$).

If $F\in k[C]$,  $L(F)$ will denote the linear component of $F$, i.e. the sum of monomials of (usual) degree 1 that appear in $F$. In the same way if $\Aid\subset k[C]$ is an ideal $L(\Aid)=\langle L(F) /\ F\in \Aid\rangle$ denotes the $k$-vector space of the linear components of elements in $\Aid$.  Obviously, if $F$ and  $\Aid$ are $\Lambda$-homogeneous, then also $L(F)$ and $L(\Aid)$ are.  If moreover $\Aid\neq k[C]$, then no element in $\Aid$ has a non-zero constant term.

 \begin{Definition}  Let $\Aid$ a $\Lambda$-homogeneous ideal in $k[C]$ and let $C'$ be a subset of $C$ and $C''=C\setminus C'$. 
 \begin{itemize} 
 \item[i)]  $\Aid/C'$ is
  the image of $\Aid$ in the quotient ring $k[C]/(C')\cong k[C\setminus C']$; we will  say that $\Aid/C'$ is obtained  \emph{specializing} $C'$ to 0.
 \item[ii)] $C'$ is a \emph{set of eliminable variables} if, for every $c'\in C'$,   $L(\Aid)$ contains elements of the type $c' +l_{c'}, $ with $f_{c'}\in k[C'']$.
   $C'$ is a \emph{maximal set of eliminable variables} if  $\vert C'\vert = \dim_k(L(\Aid))$.
 \item[iv)] for a set  $C'$  of eliminable   variables,  $ \Aid\setminus C'$ denotes the ideal $\Aid \cap k[C \setminus C'] \subset k[C\setminus C']$; we will say   that $ \Aid\setminus C'$  is obtained  \emph{eliminating} $C'$.
  \item[v)] \emph{the embedding dimension} $ed(\Aid)$ of $\Aid$ is the difference $s-\dim(L(\Aid))$ that is the cardinality of the complementary $C''=C\setminus C'$ of a maximal set of eliminable variables.
  \item[vi)] the \emph{dimension}   $\dim(\Aid)$ of $\Aid$ is the Krull dimension of the quotient ring $k[C]/\Aid$.
 \end{itemize}
\end{Definition}
If $c \in L(\Aid)\cap C$, then it belongs to every maximal set of eliminable variables. If
 $c $ does not appear in any element of $ L(\Aid)  $, then  it  belongs to the complementary of every set of eliminable variables. Finally if $c$ does not belong to $ L(\Aid)$ but   appears in some elements of it, then we can find a maximal set of eliminable variables containing $c$ and another one not containing $c$.
\begin{Proposition}\label{lambda}  Let  $\Aid \subsetneq k[C]$ be a $\Lambda$-homogeneous ideal and let $C'=\{c'_1, \dots, c'_r\}$ a set of eliminable variables for $\Aid$ and $C''=C\setminus C'$. 
 Then $\Aid$ has a set of generators of the following form:
\begin{equation}\label{generatori}\Aid=( c'_1+g_1, \dots, c'_r+g_r, f_{1}, \dots, f_d) \end{equation}
 with  $g_i, f_j \in k[C'']$, 
 and there is a natural isomorphism: 
\begin{equation}\label{quoziente} k[C]/\Aid \cong k[C'']/(\Aid \cap  k[C''])= k[C'']/(f_{1}, \dots, f_d)  .\end{equation}
If    moreover $r=dim(L(\Aid))$, then $f_j\in (C'')^2 $.

A set of generators  (\ref{generatori})   can be obtained as the reduced \Gr basis of $\Aid$ with respect to any elimination order of the variables $C'$.

\end{Proposition}

\begin{proof}  
We can  assume, up to a permutation of indexes, that $C'=\{c_1, \dots, c_r\}$ and that $\Lambda(c_1) \preceq \Lambda(c_2)\preceq  \dots \preceq \Lambda(c_r)$.   Now we proceed recursively. First of all let us   choose  a homogeneous polynomial $c_1+h_1, \in \Aid$ whose  linear forms  are $c_1+l_1\in L(\Aid)$ with $l_1\in k[C'']$ and observe that     
 $h_1\in k[C'']$: in fact neither $c_1$, nor the other variables $c_i\in C'$ can divide monomials of degree $>1$ that appear in $h_1$, because all the monomials in $h_1$ have the same $\Lambda$-degree.  So, we set $g_1=h_1$ and apply the same argument to the ideal $\Aid \cap k[C\setminus \{c_1\}]$ in $k[C\setminus \{c_1\}]$. In this way we obtain the $c_i+g_i$'s. Note that $k[C]/\Aid \cong k[C\setminus \{c_1\}]/\Aid \cap k[C\setminus \{c_1\}]$.

We can complete them  to a $\Lambda$-homogeneous set of generators  for $\Aid$:  $c_1+g_1, \dots, c_r+g_r, h_1, \dots, h_d$ and eliminate $c_1, \dots, c_s$ in the polynomials $h_j$ obtaining the $f_1, \dots, f_d$.

Finally  let $\prec'$ an elimination order of the variables $c_1, \dots, c_r$ in $\Tc$. It is quite evident that  Buchberger  algorithm applied to the set of generators given in (\ref{generatori})   only  changes (if necessary)  $f_{1}, \dots, f_d$ with a \Gr basis   of the ideal they generate $\Aid \cap k[C'']$ of $k[C'']$.
\end{proof}

In the following it will be useful the following technical  result. 

\begin{Lemma}\label{pereliminare} Let  $\Aid \subset k[C]$ be a $\Lambda$-homogeneous ideal and  $C'$
a subset of $C$ with cardinality   $s$. Then:
\begin{itemize}
\item[i)]     $\dim(\Aid)-s\leq \dim(\Aid/C')\leq \dim(\Aid)$ and $\ed(\Aid)-s\leq \ed(\Aid/C')\leq \ed(\Aid)$. 
\item[ii)] If $C' $ is a set of eliminable variables, then $\dim(\Aid \setminus C')= \dim(\Aid)$ and   $\ed(\Aid\setminus C')=\ed(\Aid)$.
\item[iii)]   If  no element of $C' $ appears in $L( \Aid)$, then  $\ed(\Aid/C')= \ed(\Aid)-s$.
\end{itemize}\end{Lemma}

\begin{proof} The statements follow from   well known facts in linear and commutative algebra  and   from what proved in Proposition \ref{lambda}.  
\end{proof}

We can now apply the above results to the ideal $\Aid$ of the stratum $\St(J)$.

\begin{Corollary}\label{usolivelli} Either the stratum $\St(J)$ of a monomial ideal $J$   is (isomorphic to) an affine space  or the origin, that is the point corresponding to $J$,  is  singular.

Moreover $\St(J)\cong \Af^d$   if and only if $d=\ed(\Aid(J))=\dim(\Aid(J))$.
\end{Corollary}
\begin{proof} If $C''$  is a maximal set of $s$ eliminable  variables for $\Aid(J)$, then $\St(J)$ can also be seen as the closed subscheme in $\Af^d $ defined by the ideal $\tilde{\Aid}(J)=\Aid(J)\cap k[C \setminus C']$, where $d=\vert C\vert -s$. If $\tilde{\Aid}(J)=(0)$, then $\St(J)$ is the affine space $\Af^d $, while if $\tilde{\Aid}(J)\neq (0)$, then the origin is a singular point for $\St(J)$ as a subscheme of $\Af^d $, because the tangent cone  is not contained in any hyperplane.
\end{proof}

Using the  properties that we have just proved, we are now able to sketch an effective and  easily implementable procedure to obtain a set of generators for the ideal $\Aid(J /  \tau)$ that defines   the stratum $\St(J /  \tau)$ of a monomial ideal $J\subset k[X]$ with respect to a reliable term order $\prec$.

\begin{procedure}\label{procedura} Let $J$ be a monomial ideal in $k[X]$ with monomial basis $B_J=\{X^{\alpha_1}, \dots, X^{\alpha_r}\}$ and let $\tau$ be a list of subsets of $k[X]$.
\vspace{-3pt}
\begin{itemize}
	\item  construct the polynomials $\{F_1, \dots, F_s\}$ by appending to each monomial $X^{\alpha_i}$ in $B_J$   a linear combination of the monomials in $\tail_{\tau\cup J}(X^{\alpha_i})$  where  the coefficients $C$ are variables;
\item for each pair of indexes $i,j$ consider the $S$-polynomial $S_{ij}$ of $F_i$ and $F_j$ and   reduce it completely to $\tilde{H}_{ij}$ w.r.t. $B_J$; 
\item write the polynomials  $\tilde{H}_{ij}$'s collecting monomials in the variables $X$ and extract  their coefficients in $k[C]$: they are      a set of generators for  $L(\Aid(J /  \tau))$;
\item  choose a  set  $C'$ of $\dim(L(\Aid(J /  \tau)))$  eliminable parameters (i.e. a set of generators for the initial ideal of $L(\Aid(J /  \tau))$ with respect to any term ordering); 
\item  reduce $S_{ij}$  to $H_{ij}$ w.r.t. $\{F_1, \dots, F_s\}$ and to a term order in $\Txc$ which is an elimination order of the variables  $X$ and is given by $\prec$ on $\TT_X$;
\item write the polynomials  $H_{ij}$ collecting monomials in the variables $X$ and extract  their coefficients in $k[C]$: they generate the ideal $\Aid(J /  \tau)$.
\item Finally: in order to obtain an  ideal defining $\St(J /  \tau)$ as an  subscheme of the  affine space of minimal dimension $\Af^{N}$ where $N=\ed(\Aid(J /  \tau))$, compute $\Aid(J /  \tau)\cap k[C\setminus C']$  using an elimination term order of the $C'$.
\end{itemize}
\end{procedure}

\begin{remark}\label{sizigie} Thanks to a well known property of \Gr bases, we can obtain the  ideal of $\St(J /  \tau)$ only  considering  the $S$-polynomials $H_{ij}$ corresponding to a set of generators for the first syzygies of  $B_J$ (see for instance \cite{CLO}, Theorem 9, page 104).
\end{remark}
 Even in quite simple cases, the number of parameters $C$ that appear in the $F_i$    is in general very high and furthermore  the ideal $\Aid(J)$ requires a big amount of   generators. For these reasons   it is not simple to compute a \Gr basis for $\Aid (J)$ with respect to any  term order. The main advantage of our procedure is that it enables us to know in advance a maximal set of eliminable parameters, only requiring some linear algebra  computations so that we can substantially reduce the number of variables involved in the computations.

In fact, if $f\in k[C]$ is the coefficient in $H_{ij}$ of a monomial $X^\beta \in\Tx$, then every monomial in $f$  has level  $\lambda=\Lambda(\lcm(X^{\alpha_i},X^{\alpha_j}))-\Lambda(X^\beta)$ and every coefficient $c_{i\beta}$ that appears in it has level lower or equal to $\lambda$.  Then we can obtain a set of generators for $\Aid (J)$ \lq\lq level by level\rq\rq (starting from the lowest one),  after  specialization to 0 of every parameter of level higher than $\lambda$.

\medskip

In the following sections  we will see both  explicit examples and  theoretical applications of  Procedure  \ref{procedura}.

\section{Examples}\label{sec:esempi}

The explicit computations of strata presented in the following examples have been performed using the Maple12 procedure  available at:

\noindent{\texttt{http://www2.dm.unito.it/paginepersonali/roggero/InitialIdeal(Maple)/}}

 As already observed, we can think of $\St (J/ \tau)$ as a  section of  $\St (J)$ with a linear subspace. The following example shows that  $\St (J /  \tau)$ can   be singular at the origin, even if $\St (J)$ is isomorphic to an affine space and vice versa. 

\begin{Example}\label{tagliare} Let us consider the ideals $J=(x^2 y, xy^2)$, $  J'=(y^3,z^3) \subset k[x,y,z,t]$,  and     the term order \emph{\texttt{RevLex}} with $x\succ y\succ z\succ t$. Explicit computations show that    $\St(J)$ of $J$  is smooth and more precisely is isomorphic  to the affine space $\Af^{12}$. If we cut $\St(J)$  with the codimension $4$ linear space   given by $c_{1,y^3}=c_{1,z^3}=c_{2,y^3}=c_{2,z^3}=0 $,  we obtain the stratum $\St(J /  J')$. Again by an explicit computation one can see that  $\St(J /  J')$ has two irreducible components, one of dimension $10$ and one of dimension $9$. The elimination of all the eliminable variables gives an embedding of $\St(J /  J')\subset \Af^{11}$  such that  the two components become respectively a quadric hypersurface and  a codimension 2 linear  space, both through the origin.
\end{Example}

\begin{Example} Let $J=(x^2 ,xy, xz,y^2)$ be the ideal in $ k[x,y,z,t]$ with the term order \texttt{Revlex} such that $x\succ y\succ z\succ t$. The homogeneous stratum $\St_h(J)$ is a reducible variety with a component of dimension 11 and a component of dimension 8 both containing  the origin.  If  $J'=(t)$, we obtain $\St_h(J / J')\cong \Af^4$.  

Observe that in a natural way $\St_h(J)\cong \St(J_1)$ and  $\St_h(J /  J')\cong \St_h(J_1)$ where $J_1=(x^2 ,xy, xz,y^2)$ is an ideal in $  k[x,y,z]$. Then this same example also shows that the homogeneous stratum can be smooth, while the stratum is not.
\end{Example}

The last example shows that in general, but not always, the homogeneous stratum $\St(J,\prec)$ depends   on the ideal $J$ itself and not only on the projective scheme it defines.
\begin{Example}\label{3punti} Let us consider  in the polynomial ring $ R= k[x,y,z,t]$ with $x\succ y\succ z\succ t$ the ideal: $$J=(x^3 ,x^2y, x^2z,x^2t,xy^2,xyz,xyt,xz^2,xzt,xt^2,y^3,y^2z,y^2t,yz^2,yzt,yt^2,z^3)$$ 
and its saturation $J_{sat}=(x,y,z^3)$. Of course 
   the projective schemes $\Proj(R/J)$ and $\Proj(R/J_{sat})$  are canonically isomorphic, while  the homogeneous strata   of the two ideals may not  be. In fact,  if we fix the   term order  \texttt{Lex} and compute the  strata, we obtain  $\St_h(J,  \texttt{Lex})\cong\Af^9$, while $\St_h(J_{sat},  \texttt{Lex})\cong \Af^7$.  
   
 On the other hand,  if we choose     the term  order   $\texttt{RevLex}$,   the homogeneous strata of $J$ and of $J_{sat}$  are both isomorphic to  $ \Af^7$.
\end{Example}

\section{Strata of    \texttt{Lex}-segment ideals}\label{sec:lex}

In this section we consider any  homogeneous  saturated  segment ideal with respect to the lexicographic term order $\texttt{Lex}$. Now $R=k[x_0, \dots, x_n]$ and we fix the order on the variables $x_0 \succ  x_1 \succ  \dots \succ x_n$. Every saturated  monomial ideal $J$ which is a segment  with respect to $ \texttt{Lex}$ is completely determined by the monomial in its monomial basis which is minimal with respect to $\prec_{\texttt{Lex}}$: if  $\mathbf{m}=x_{0}^{a_{n}}x_{1}^{a_{n-1}}\cdots x_{n-2}^{a_{2}}x_{n-1}^{a_{1}}$ is such a monomial, and $r=\sum a_i$,   then $J_r$ is the vector space generated by all the monomials in degree $s$ that are $ \texttt{Lex}$-higher than 
$\mathbf{m}$ and $J$ is the saturation of the ideal generated by $J_r$.  Following \cite{RS} we then identify every such monomial with the list of exponents $(a_n, \dots, a_1)$ (note that our use of indexes is slightly different from that of  \cite{RS}).  More generally:

\begin{Definition}\label{def:idealeJ}  Let us consider an integer  $q \leq n$   and  any sequence of  $q$  non negative integers $ (a_q,...,a_{1 })$  with $\sum a_i\neq 0$.  Let  $s\geq 0$ be  the minimal index such that $a_{s+1}\neq 0$. We will denote    by $J(a_q, \dots, a_1) $ the ideal in $R$ generated by the set of $q-s$ monomials $\m_1, \dots, \m_{q-s}$ where:
\begin{equation}
  \m_i= x_{n-i}\cdot \prod_{j=i}^q x_{n-j}^{a_j}  \quad  \hbox{ if }  \quad  1 \leq i <q-s \quad  \hbox{ and } \quad \mathbf{m}_{q-s}=  \prod_{j=s+1}^{q} x_{n-j}^{a_j} 
\end{equation}
\end{Definition}

\begin{Example}\label{casibassi} For $q=n$,  $J(0,\dots, 0, a_1)$ is the ideal $(x_1,x_2, \dots, x_{n-2}, x_{n-1}^{a_1})$  defining a non-reduced structure of length $a_1$ over the origin in $\PP^n$.  For a general $q$, $J(0,\dots, 0, a_1)=(x_{n-q}, \dots, x_{n-2},x_{n-1}^{a_1}) $ defines   a non reduced structure over a linear space of codimension $q$. 

For every $q> 1$,  $J(a_q, 0,\dots, 0)=(x_{n-q}^{a_q})$  defines the hypersurface in $\PP^n$ corresponding to the divisor $a_{q}H$ where $H$ is the hyperplane $x_{n-q}=0$. 
\end{Example}
      In \cite{RS}  it is proved that very    saturated lex-segment ideal $J$ of $R$ is  of the type
$ J(a_n,...,a_1)$ and that it   defines a  subschemes of $\PP^n$ which   corresponds to   smooth points  in the Hilbert scheme $\hilbp$, where $p(z)$ is the Hilbert polynomial of $R/J$. By the universal property of Hilbert schemes, $\St_h(J)$ can be embedded as a locally closed subscheme in $\hilbp$. In \cite{NS} it is claimed  that the smoothness of   $\St_h(J)$ can be deduced   from that of the Hilbert scheme containing it. This argument clearly holds  when   $\St_h(J)$ corresponds to an open subset of $\hilbp$. However, it is not difficult to verify that $\St_h(J)$ has in general a strictly  lower dimension than the component of the Hilbert scheme containing it, hence it cannot correspond to an open subset of  $\hilbp$ (see Remark \ref{ultima} or,  for an explicit computation in a special case, Example \ref{3punti}).

\medskip

 The two ideals $J= J(a_n,...,a_1)$ and $ I=J(a_n,...,a_1)_{\geq r}$ (where $r=\sum a_i$) define the same point in the corresponding Hilbert scheme. In  \cite{LR} it is proved that the homogeneous stratum $ \St_h(I, \texttt{Lex})$  is scheme-theoretically isomorphic to an open subset of the  Hilbert scheme and that it is also  isomorphic to an affine space $\Af^d$ of a suitable dimension $d=d(a_n,...,a_1)$.  On the other hand, the homogeneous stratum  $ \St_h(J, \prec)$, with respect to any term ordering $\prec$, is the locally closed subscheme of $ \St_h(I, \texttt{Lex})$  cut by a   suitable linear space,   because it can be realized as a stratum of the type $\St_h(I /  \tau)$.  Even if $\St_h(I,\prec)$ is smooth and so isomorphic to an affine space,  nothing can be said in general about the smoothness  of a stratum of the type $\St_h(I /  \tau)$, because it can be singular or even reducible and non-reduced  (see   Example \ref{escluso}). 
 
 \medskip

 However, in the case of the ideals $ J(a_q,...,a_1)$ we will prove that  their strata w.r.t. any reliable term ordering and their homogeneous strata w.r.t. any term ordering are always isomorphic to   affine spaces. As a consequence this holds for the  strata of any lexicographic saturated ideal. If the fixed term ordering is \texttt{Lex}, we also give a formula for the dimension of $\St_h( J(a_n,...,a_1), \texttt{Lex})$.

The starting point of our proof is    the  nice \lq\lq inductive structure\rq\rq\ of the ideals of the type $J(a_q,...,a_1)$, underlined in \cite{RS}:
\begin{equation}\label{ind1}  J(a_q,...,a_{1})=x_{n-q}^{a_q}\cdot((x_{n-q})+J(a_{q-1},...,a_{1}))\end{equation}
More generally, we will   consider  ideals having  a  structure of the type (\ref{ind1}) and study their strata. For the meaning of $\St(J_0/(\n_0))$ we refer to Example \ref{escluso}.

\begin{Theorem}\label{torno} Let $X, Y $ be two disjoint sets of variables and let $\prec$ be    any reliable  term order  in $\Txy$.  
Let us consider in $k[X,Y]$ a  monomial ideal $J_0$ with basis  $\n_1, \dots , \n_r\in \Tx$,  two monomials  $\m,\n_0\in \Ty$ with  and  the monomial ideals  $J_1=((\n_0)+J_1)$. If $\St(J_0/(\n_0)) \cong \Af^{N_0}$ and either $\n_0$  has degree $1$ or $\prec$ is an  elimination order  of the set of variables $Y$,
then the stratum   of the ideal:
 $$J= \m((\n_0)+J_0)$$
is isomorphic to an affine space. More precisely:
$$\St(J) \cong \Af^{N_0+N_1+N_2}  $$
where $N_1=\vert \tail(\m)\vert$ and $N_2=\vert \tail(\n_0)\setminus J_0\vert$. 

The same holds for $\St_h(J)$ assuming that  $\St_h(J_0/(\n_0)) \cong \Af^{N_0}$,  $N_1=\vert \htail(\m)\vert$,  $N_2=\vert \htail(\n_0)\setminus J_0\vert$.
\end{Theorem}

\begin{proof} We prove the statement only for what concerns the stratum, the homogeneous case being analogous.

In order to compute the stratum of $J_1$ (and respectively the stratum of $J$)  we follow   Procedure \ref{procedura}   starting from   polynomials $G_0,G_1  \dots, G_r$ (respectively $F_0, F_1, \dots, F_r$) as in (\ref{polF}).  If:
  \begin{eqnarray*} G_0=\n_0+\sum D_{s}\vv_{0s}  &  \\
  G_i=\n_i+\sum C_{is}\vv_{is}  &    \hbox{ if } \quad i\geq 1
  \end{eqnarray*}
  where $\vv_{is}\in \tail_{J_1} (\n_i)$, then we can write the  $F_i$    collecting  first   all the monomials multiple of $\m$ and, among the remaining ones, all the monomials multiple of $\n_i$, in the following way:
 \begin{eqnarray*}
 F_0=\m G_0+ \n_0\left(\sum C'_{0h}\p_h\right)+\sum C''_{0t}\q_{0t} & \\
 F_i=\m G_i+ \n_i\left(\sum C'_{ih}\p_h\right)+\sum C''_{it}\q_{it} \ &     \hbox{ if } \quad i\geq 1
 \end{eqnarray*} 
 where    $\p_h\in \tail(\m)$ (if $i=0$ also $\n_0 \p_h \notin (\m)$) and $\q_{it} \in \tail(\m\n_i)\setminus (\n_i,\m)$. 
 
 First of all we observe that there is a $1-1$ set-theoretic correspondence between  $\St(J)$ and the product $\St(J_0 / \n_0) \times \St((\m))\times \St((\n_0)/J_0)\cong \Af^{N_0} \times \Af^{N_1} \times \Af^{N_2}$.   If $I\in \St(J)$, that is if $\In(I)=J$, then $I$ and $J$ share the same Hilbert polynomial and then $\htt(I)=\htt(J)=1$. If $I=F\cdot I_1$ and $\htt(I_1)\geq 2$, then $\LM(F)=\m$ and $\In(I_1)=J_1$. On the other hand if $\LM(F)=\m$ and $I_1$ is an ideal such that $\In(I_1)=J_1$, then $\In(F\cdot I_1)$ has the same Hilbert function than $J$ and $\In(F\cdot I_1)\supseteq J$, hence  they coincide. Hence  the  points in $\St(J)$ correspond $1-1$ to  points in $\St(\m)\times \St(J_1)\cong \Af^{N_1} \times \St(J_1)$. By   hypothesis, $\n_0$ is prime to each $\n_i$ so that the reductions of the $S$-polynomials $S(G_0,G_i)$  w.r.t.   $G_0, \dots, G_r$ vanish and so they do not contribute to the construction of $\Aid(J_1)$. Moreover our assumptions on   $\n_0=y$ insure that  it cannot appear in any step of reduction of $S(G_i,G_j)$ for every $i,j>1$. Hence   $ \St(J_1)\cong \St((\n_0)/  J_0)\times  \St(J_0 / (\n_0)) \cong \Af^{N_2+N_0}$.

  Now we shall prove that this correspondence is in fact a scheme-theoretical isomorphism, proving that the embedding dimension of $\St(J)$ is lower than or equal to $N_0+N_1+N_2$ computing a suitable set of elements  of $L(\Aid(J))$ (for the meaning of $L(\Aid(J))$ see \S \ref{sec4}).  More precisely we will prove that  all the parameters are eliminable for $\Aid(J)$ except:
\begin{description}
\item[- ] $N_0$ among the $C_{is}$ and the $N_2$ parameters $D$;
\item[- ]   $N_1$ parameters $C'$, for instance the $C'_{rh}$.
\end{description}

The first step is to prove    that $\Aid(J_0/(\n_0))\subseteq \Aid(J) $, where 
 $\Aid(J_0/(\n_0)) $ is the ideal in $k[C]$ that we obtain 
applying   Procedure \ref{procedura} to the polynomials $G_1, \dots, G_r$.  We can write $S(F_i,F_j)$ as $ \m S(G_i,G_j)+M_{ij}$ where  $M_{ij} $ does not contain monomials multiple of $\m$. By hypothesis, $\m \n_0$ does not divide any monomial in $\m S(G_i,G_j)$ and no monomial in $M_{ij}$ belongs to $J$. Hence the reduction of $S(F_i,F_j)$ with respect to $J$ can be written as $H(F_i,F_j)=\m H(G_i,G_j)+\tilde{M}_{ij}+M_{ij}$, where $H(G_i,G_j)$ is the reduction of $S(G_i,G_j)$ with respect to $G_1, \dots G_r$ and no monomial in $\tilde{M}_{ij}$ is multiple of $\m$.  Then, collecting the coefficients of monomials that are multiple of $\m$ we obtain the wanted inclusion  
  $\Aid(J_0/(\n_0)) \subset \Aid(J)$. 
  
 \medskip
 Now we fix a  maximal set of eliminable variables for $\Aid(J_0/(\n_0))$ and specialize to $0$ the complementary set of $N_0$ variables
 $C$. Thanks to what just proved, this is equivalent to specialize to $0$ all the variables $C$. By  Lemma \ref{pereliminare} we then have  $ \ed(\Aid(J))  \leq \ed(\Aid(J)/C)+N_0$. 
 
 As $S(F_0,F_i)=\n_iF_0-\n_0F_i$ and $\m\n_i\vv_{0s}\in J$, no variable $D$ appears in $L(\Aid(J))$ and so we can  specialize to $0$ also the $N_2$ variables $D$ obtaining that $ \ed(\Aid(J))\leq \ed(\Aid(J)/C\cup D)+N_0+N_1$ and that the procedure can be continued using the polynomials:   
 \begin{eqnarray*}
 \tilde{F}_0=\m\n_0+ \n_0\left(\sum C'_{0h}\p_h\right)+\sum C''_{0t}\q_{0t} & \\
 \tilde{F}_i=\m \n_i+ \n_i\left(\sum C'_{ih}\p_h\right)+\sum C''_{it}\q_{it} \ \ &     \hbox{ if } \quad i\geq 1
 \end{eqnarray*} 
Let us  observe that  two  parameters of the same level are both in $C'$ or both in $C''$, because  the hypothesis on the $\q_{jt}$  does  not allow equality between $\m \n_i/\q_{it}$ and    $\m \n_{j}/\n_j\p_h$. Moreover if  $j>0$ there are no parameters in $\tilde{F}_j$ of  the same level of some $C''_{0t}$.  
 
Computing the reduction of $S(\tilde{F}_0,\tilde{F}_r)$ with respect to $J$, we can see that for every $t$ the  parameter   $C''_{0t}$ belongs to $L(\Aid(J))$ and for every $h$ either   $C'_{0h}$ or the difference $C'_{0h}-C'_{rh}$ belongs to $L(\Aid(J))$: in fact    $\n_0\n_r\p_h \notin J$ and $ \n_r \q_{0t} \notin J$.   Moreover using the $S$-polynomials $S(\tilde{F}_i,\tilde{F}_r)$ with $i>0$ we see that also $C'_{ih}-C'_{rh} \in L(\Aid(J))$.

Finally, we use the $S$-polynomials $S(\tilde{F}_i,\tilde{F}_j)$ with $i,j>0$ to prove that also the other  parameters $ C''$ belong to $L(\Aid(J))$. We may assume by simplicity that $i=1$.    If  there is an index $j>0$ such that
no parameter of the same level as $C''_{1t_1}$ appears in   $\tilde{F}_j$, we can see    that $ C''_{1t_1}\in L(\Aid(J))$ using  the reduction of $S(\tilde{F}_1,\tilde{F}_j)$  with respect to $J$.
 If, on the contrary, for every $j>0$ there is a parameter $C''_{jt_j}$ of the same level as $C''_{1t_1}$  (so that $\n_1\q_{jt_j}=\n_j\q_{1t_1}$), then using $S(\tilde{F}_1,\tilde{F}_j)$  we see that $C''_{1t_1}-C''_{jt_j}\in L(\Aid(J))$ for every $ j>0$. Moreover there is at least an index  $j_0>0$ such that $\n_0\q_{j_0t_{j_0}}  \notin J$; in fact, $\n_0\q_{jt_j} $ is not a multiple of $\n_0 \m$ and the inclusion $\n_0(\q_{1t_1}, \dots, \q_{rt_r})\subseteq \m(\n_1, \dots,\n_r)$ would imply $\n_0\n_1(\q_{1t_1}, \dots, \q_{rt_r})=\n_0\q_{1t_1}(\n_1, \dots,\n_r)\subseteq \m\n_1(\n_1, \dots,\n_r)$ and then $\q_{1t_1} \in  (\n_1)$ against the hypothesis on $\q_{1t_1}$. So, using $S(\tilde{F}_0,\tilde{F}_{j_0})$, we see that  $C''_{j_0t_{j_0}}\in L(\Aid(J))$, hence also $C''_{1t_1}\in L(\Aid(J))$. 
 
\medskip

Thus $\ed(\Aid(J))\leq N_0 +N_1+N_2$ and on the other hand $\ed(\Aid(J))\geq \dim(\Aid(J)) =N_0 +N_1+N_2$ and this allows us to conclude.
\end{proof}

The following example shows that Theorem \ref{torno} does not hold if we weaken the hypotheses on both $\n_0$ and $\St((\n_0)/J_0)$.  

\begin{Example}
Let us consider the ideal $J=(y^2, x_1^2,  x_2x_3, x_1x_2^2, x_1x_3^2)$ in $k[y,x_1,x_2,x_3]$ and the term ordering $\texttt{Lex}$ with $x_1\succ y\succ x_2\succ x_3$. This ideal is of the type $ \m((\n_0)+J_0)$ with  $\m=1$, $\n_0=y^2$ and $J_0=(x_1^2, x_1x_2^2, x_1x_3^2, x_2x_3)$. Then $N_1=0$ and  $N_2=\vert \htail(\n_0)\setminus J_0\vert =4$. By a direct computation we can prove  that $\St_h(J_0/(\n_0))$ is isomorphic to the union of $4$ distinct linear spaces of dimension $6$ and also that $\St_h(J)$ is isomorphic to $ \Af^{12}$. Then $\St_h(J)$ is not isomorphic to   $\St_h(J_0/(\n_0))\times \Af^4$.
\end{Example}
As a consequence of the previous general result we can prove that the strata and the homogeneous strata of ideals of the type $J(a_q, \dots, a_1)$, and so especially those of the $\texttt{Lex}$ segment ideals, are isomorphic to affine spaces. We point out that   Theorem \ref{torno} does not require any special condition on the order induced by the chosen term ordering on the set of variables. In the special case of the   $\texttt{Lex}$ segments ideals, if we fix the term ordering $\texttt{Lex}$ that orders in the usual way the set of variables, we can also give a formula for the dimensions of the homogeneous strata.

\begin{Corollary}\label{segmento1} Let us fix in $k[x_0, \dots, x_n]$ any  term order and any sequence of $q\leq n$ non-negative integers  $a_q, \dots, a_1$. Then    the homogeneous stratum of $J(a_q,...,a_{1 })$ is isomorphic to an affine space. The same holds for the stratum $\St(J(a_q,...,a_{1 }))$ if it is defined, i.e. if the term ordering is reliable.
\end{Corollary}
\begin{proof} We prove the statement by induction on $q$ for the homogeneous stratum, the non homogeneous case being analogous. If $q=1$ the ideal $J(a_1)$ is $(x_{n-1}^{a_1})$ and    its homogeneous stratum is an  affine space.

Let us assume that the claim holds for every $q'<q$, ($q\geq 2$) and let us apply Theorem \ref{torno} to $J_0=J(a_{q-1},\dots, a_{1})$, $\n_0=x_{n-q}$ and $\m=x_{n-q}^{a_q}$. 

Observe that in this case $\St_h(J(a_{q-1}, \dots,a_1),(x_q))$ is  isomorphic to  the stratum of $J(a_{q-1},\dots, a_{1})$  as an ideal of $k[x_0, \dots, x_{n-q-1},x_{n-q+1}, \dots, x_n]$, and so by  inductive hypothesis  it is isomorphic to an affine space. 
\end{proof}

We stress that $J(a_n,...,a_1)$ is a $\texttt{Lex}$-segment ideal with respect to the $\texttt{Lex}$ term ordering with $x_0\succ x_1 \succ \dots \succ x_{n}$. However, our result on $ \St_h\left(J(a_n,...,a_1), \prec\right)$ and $ \St(J(a_n,...,a_1), \prec)$   holds for any term ordering, even if it changes the order on the set of variables. When the term order is \texttt{Lex} and $x_0\succ x_1 \succ \dots \succ x_{n}$, we give a formula for the dimension of the homogeneous stratum.

\begin{Corollary}\label{segmento2} Let us fix in $k[x_0, \dots, x_n]$ the term ordering $\texttt{Lex}$ with $x_0 \succ x_1 \succ \dots \succ x_n$  and any sequence of $q\leq n$ non-negative integers  $a_q, \dots, a_1$. Assume that   $\sum a_i\neq 0$   and let $s\geq 0$ be  the minimal index such that $a_{s+1}\neq 0$. 

Then    $\St_h(J(a_q, \dots, a_1),\texttt{Lex})$ is isomorphic to an affine space of dimension
   $M=M(a_q, \dots, a_1)$  given by:
    $$M(a_q, \dots, a_1)=\frac{(q-s)^2-(q-s)-2}{2}+\sum_{j=s+1}^q {a_j+j\choose j}  -\sum_{j=s+1}^{q} \nu_j $$
where,  for every $j>s+1$, $\nu_j$ is such that $a_{j-1}=\dots=a_{j-\nu_j}=0$ and $a_{j-\nu_j-1}\neq 0$.
\end{Corollary}
\begin{proof} We prove the statement by induction on $q-s$. If $q-s=1$, then    the ideal $J(a_q,0, \dots, 0)$ is $(x_{n-q}^{a_q})$ and its homogeneous stratum is an affine space whose  dimension is the number $\binom{a_q+q}{q}-1$  of the elements in  the homogeneous tail of $  x_{n-q}^{a_q} $. This number coincides with    $M(a_q, 0, \dots, 0)$.

Now let us assume   $q-s\geq 2$ and that the thesis holds for lower cases of the difference $q-s$. We can apply Theorem \ref{torno} to $J_0=J(a_{q-1},\dots, a_{1})$, $\n_0=x_{n-q}$ and $\m=x_{n-q}^{a_q}$: note that the  number $s$ is the same for the two ideals  $J(a_{q},\dots, a_{1})$  and $J(a_{q-1},\dots, a_{1})$. Now we observe that    $x_q$ cannot appear in the tail w.r.t. $\texttt{Lex}$ of any monomial in the  basis of $J(a_{q-1}, \dots,a_1)$ because $x_0 \succ x_1 \succ \dots \succ x_n$. Thus $\St_h(J(a_{q-1}, \dots,a_1),(x_q))=\St_h(J(a_{q-1},\dots, a_{1}))$.   The difference between  $\dim(\St_h(J(a_{q},\dots, a_{1})))$ and  $\dim(\St_h(J(a_{q-1},\dots, a_{1})))$  is given by the sum of two numbers. The first one is the number of elements in the tail of $\m=x_{n-q}^{a_q}$, that is ${a_q+q\choose q}-1$. The second one is    the number of monomials in the tail of $x_q$   non contained in $J(a_{q-1}, \dots,a_1)$, that is   $q-\nu_q$. This is precisely the difference
 $M(a_{q}, \dots,a_1)-M(a_{q-1}, \dots,a_1)$. 
\end{proof}

\begin{remark}\label{ultima} For every    admissible Hilbert polynomial $p(z)$ in $\PP^n$,  the Hilbert scheme $\hilbp$  contains  a point, usually called \texttt{Lex}-point, corresponding to $\Proj(k[x_0, \dots, x_n]/J(a_n,\dots, a_1))$ for a suitable  sequence   $(a_n,\dots, a_1)$.     In \cite{RS} it is proved that   the  \texttt{Lex}-point is a smooth point of $\hilbp$ and an explicit formula is given for the dimension of the  \texttt{Lex}-component (also called the Reeves and Stillman component), that is the only component of   $\hilbp$  containing the \texttt{Lex}-point.     A direct comparison clearly shows that this dimension and that of $\St_h(J(a_n,\dots, a_2,0), \texttt{Lex})$ given in the previous result are in general different. 

As an evidence we can consider   $J(0, \dots, 0,a_1)=(x_0, \dots, x_{n-2}, x_{n-1}^{a_1})$ with constant Hilbert polynomial $p(z)=a_1$: the dimension of the homogeneous stratum  $\St_h(J(0,\dots, 0,a_1), \texttt{Lex})$ is $2n-2+a_1$, while the dimension of  the \texttt{Lex}-component of   $\hilb^n_{a_1}$      is $na_0$ (see \cite{RS}, especially  the remark after  Lemma 4.2).  Of course, the ideals corresponding to points of   $\St_h(J(a_1, \dots, a_r), \texttt{Lex})$ also correspond to points in   the  \texttt{Lex}-component of $\hilbp$. Then, if the dimensions are different, the dimension of  $\St_h(J(0,\dots, 0,a_1), \texttt{Lex})$   must be lower than that of the \texttt{Lex}-component. This means that in the stratification of the Hilbert scheme through strata given in \cite{NS}, the stratum of the lexicographic point is in general a   \lq\lq slice\rq\rq\ in its component, that is a locally closed subscheme of it  which is not in general an open subset.
The following result describes   an important  case in which the two dimensions agree.  
\end{remark}

\begin{Corollary} Let us fix in $k[x_0, \dots, x_n]$ the term ordering $\texttt{Lex}$ with $x_0 \succ x_1 \succ \dots \succ x_n$  and consider any ideal of the type  $J(a_n, \dots, a_2,0)$.

Then   $\St_h(J(a_n,\dots, a_2,0), \texttt{Lex})$ is isomorphic to  an open subset of the \texttt{Lex}-component of the corresponding Hilbert scheme. 
\end{Corollary}
\begin{proof}
  Every ideal $J(a_n,\dots, a_1)$ is saturated  and strongly stable and  so   the lower variable $x_n$ does not appear in any  monomial of  its monomial basis. If moreover, as in the present case,   $a_1=0$ then also $x_{n-1}$ does not appear in those monomials. We can then apply  \cite[Theorem 4.7]{LR} obtaining that     $\St(J(a_n,\dots, a_2,0), \texttt{Lex})\cong \St_h(J(a_n,\dots, a_2,0)_{\geq r},, \texttt{Lex})$, where $r$ is the Gotzmann number of $p(z)$, that is the regularity of $J(a_n,\dots, a_2,0)$.  Moreover $J(a_n,\dots, a_2,0)\geq r$ is generated b the maximal monomials of degree $r$ w.r.t. \texttt{Lex} and so, by  \cite[Corollary 6.9]{LR},  $\St_hJ(a_n,\dots, a_2,0)_{\geq r},, \texttt{Lex})$ is isomorphic to an open subset of the Hilbert scheme $\hilbp$.    
\end{proof}

\section{Strata of   \texttt{RevLex}-segment ideals}\label{sec:revlex}

Let $Z$ be a general set of $\mu$ points in $\PP^2$ in general position,  $I_Z \subseteq k[x,y,z]$ be its ideal.   In \cite{CS} it is proved that the initial ideal of $I_Z$  with respect to the  the term order $\mathtt{RevLex}$ with $x\succ y \succ z$ is the ideal $\idR(\mu)=( x^r, \dots  x^{r-i}y^{i+\delta_i}, \dots y^{r+\delta_r})$, for some positive integers    $r$ and $t$  and  $\delta_i=0 $ if  $i<t$ and   $\delta_i=1 $ otherwise.  In every degree $s$ the  monomials contained in
 $\idR(\mu)_s$ are the maximal ones w.r.t. \texttt{RevLex}, that is they form a  \texttt{RevLex}-segment.  Especially, in degree $\mu$, the ideal $\idR(\mu)$ contains  the   $\binom{\mu +2}{2}-\mu$ maximal  monomials w.r.t. \texttt{RevLex}  and for this reason, using the terminology introduced in   \cite{CLMR}, $\idR(\mu)$ is a $Hilb$\emph{-segment ideal} and corresponds to the \texttt{RevLex}\emph{-point} of $\hilbp$ where $n=2$ and $p(z)$ is the constant polynomial $\mu$ (see for instance \cite{MS}). 
 
 In \cite{CLMR} this construction is generalized to the case of $\mu$ points in $\PP^n$, $n\geq 3$ proving that  $\hilb^n_N$  contains a \texttt{RevLex}-point and that it is    a  singular  point  if $\mu>n$   (see \cite[Theorem 6.4]{CLMR}). 
 
In this last section we consider another type of generalization of $\idR(\mu)$, that is    the ideals $ \idR(\mu,n)$    that are  generated by the same monomials    $ x^r, \dots  x^{r-i}y^{i+\delta_i},$ $ \dots,  y^{r+\delta_r}$  as $\idR(\mu)$   but in $k[x, y,z_1,\dots, z_{n-1}]$.  Using   Procedure \ref{procedura}   we are able to prove that the homogeneous  stratum  $\St_h(\idR(\mu,n), \texttt{RevLex})$ is isomorphic to an affine space  and  to give a formula for its dimension.  We underline that even though the case $n=2$ of   our statement is a well known fact,   we  need to prove it explicitly because the proof itself, not only the statement, is the starting point for the proof of the general case. Finally, we note that the ideals $\idR(\mu,n)$ define arithmetically Cohen-Macaulay subschemes in $\PP^n$ and so the image of their strata in the corresponding Hilbert schemes are contained in the smooth open subset studied in \cite{Elling}.

\begin{Theorem}\label{punti} 
Let $\prec$ be the    term order \texttt{RevLex} in $k[x, y,z_1,\dots, z_{n-1}]$ with $x\succ y\succ z_1 \succ \dots \succ z_{n-1}$. For any positive integer $\mu$ let     $r$ and $t$ be  the only integers such that $1\leq t\leq r+1$ and $\mu=\frac{(r+2)(r+1)}{2}-t$ and    let us consider  the ideal in $k[x, y,z_1,\dots, z_{n-1}]$: 
 $$\idR=\idR(\mu,n)=(x^r, \dots  x^{r-i}y^{i+\delta_i}, \dots y^{r+\delta_r})$$ 
 where $\delta_i=0 $ if  $i<t$ and   $\delta_i=1 $  otherwise. Then the homogeneous stratum $\St_h (\idR(\mu,n), \texttt{RevLex})$ is isomorphic to an affine space $\Af^{N}$ of dimension
\begin{equation}\label{primomodo}  N=N(n,\mu)=2(n-1)\mu+t(r+1-t)\binom{n-2}{2}.
\end{equation}
\end{Theorem}
\begin{proof} As the term ordering is fixed, in the following we  avoid to indicate  it. We first prove the assertion in a special case.

\textbf{Step 1}: the case  $n=2$.  As above said, the homogeneous stratum $\St_h(\idR)$ of $\idR=\idR(\mu)$ contains an algebraic family of dimension $2\mu$ corresponding to the sets of $\mu$ general points in $\PP^2$. Hence in order to prove that $\St_h(\idR)\cong \Af^{2\mu}$ it is sufficient to prove that $\ed(\Aid(\idR))\leq 2\mu$ (see Corollary \ref{usolivelli}).

Following  Procedure \ref{procedura}, let  $F_i=x^{r-i}y^{i+\delta_i}+\sum c_{ijk}x^jy^k z^h$, where   $j+k+h = r+\delta_i$ and   $x^jy^k \notin \idR$.
  We want to show that  all the variables $C$ are eliminable for $\Aid$ except   $2\mu$ of them and for this we compute $L(\Aid)$, where $\Aid=\Aid(\idR)$. Now we recall that we can obtain   $L(\Aid)$  only using   the $S$-polynomials $S_i=S(F_{i},F_{i+1})$,  because they correspond to a set of generators for the syzygies  of   $\idR$ (see Remark \ref{sizigie}) and their reduction $\tilde{H}_i$ w.r.t. $\idR$.   Moreover, $L(\Aid)$  has a set of generators that are linear combinations of coefficients $C$ of the same level. Now observe that any two coefficients $c_{i,j,k} \in C$ of the same level must share the same number  $h=r+\delta_i-j-k$ and that   there are no parameters of the same level in both $F_0$ and $F_r$ because no monomials $\n$ and $\m \in k[x,y]\setminus \idR$    satisfy $  \frac{x^r}{\n}=\frac{y^r}{\m}$.   
  thanks to these remarks, we can consider a fixed $h$ at a time, that is think of polynomials  $F_i$ of the type $ F_i=x^{r-i}y^{i}+z^h\sum_{j+k=r+\delta_i-h } c_{ijk}x^jy^k $ with $x^jy^k\notin \idR$.

First we take in consideration the special case  $2\mu=r(r+1)$, so that $t=r+1$ and $\idR=(x^r, x^{r-1}y,  \dots, y^r)$   and prove that  all   coefficients $C$ are eliminable, except those with $h=1$, that are $(r+1)r=2\mu$ (in this case $x^jy^k\notin \idR$ requires $h>0$). 

So, let us fix any $h\geq 2$ and let $l$ be the level of $c_{ijk}$, with $j+k=r-h$  and let   $p$ and $q$  be respectively the minimum and the maximum $i$ such that in $F_i$ there is a parameter of level $l$ that we will denote simply $\overline{c}_i$. If $p\neq 0$, using $\tilde{H}_{p-1}$, $\tilde{H}_{p}, \dots,\ \tilde{H}_{q-1}$ we find  $-\overline{c}_p$, $\overline{c}_p-\overline{c}_{p+1}, \dots, \overline{c}_{q-1}-\overline{c}_{q} \in L(A)$, so that   $\overline{c}_i\in L(\Aid)$. If $p=0$ and $q\neq r$, we argue in the same way using $\tilde{H}_{0}  \dots,\tilde{H}_{q-1}, \tilde{H}_{q}$.
This allows us to conclude.

Now we consider the   case $t\leq r$ so that $S_{i}=yF_{i}-xF_{i+1}$ if $i\neq t-1$, while $S_{t-1}=y^2F_{t-1}-xF_t$. We assume by induction that  the statement holds for $r'=r-1$, $t'=t$ and $\mu'=\frac{(r+1)(r)}{2}-t $ (the initial case $r=1$ being obvious).

 As above,  all the variables $c_{ijk}$ such that either $j+k\leq r+\delta_i-2$ or $i\neq t-1$ and $j+k=r+\delta_i-1$ are eliminable.   
  For every $i\leq r-1$ and every monomial of the type $y^{b}$,  such that $b\leq r-1$  and also  $b\neq r-1$ for $i=t-1$, using   $\tilde{H}_{i} $ we can see that $c_{i0b}\in L(\Aid)$.  Now we specialize to 0 all the parameters $c_{i0b}$ for every $i\neq r$: in this way the embedding dimension drops at most by $  r+1 $ because those coefficients  are all in $L(\Aid)$ except  at most $r+1$ of them, namely  $c_{i0r}$, $i=0, \dots, r-1$  and  $c_{(t-1) 0 (r-1)}$.

  Now  the polynomials $F_i$ take the form: $F'_0=xG_0, \dots, F'_{r-1}=xG_{r-1}$, $F'_r=xG_r$ where  $G_0, \dots, G_{r-1}$ are the polynomials that we use in order to construct the homogeneous stratum of $\idR(\mu')$. Now we observe that  $y^{r+1}$ does not divide any monomial in $S(G_{i},G_{i+1})$ if $i\leq r-2$, so that the reduction of $S(F_{i},F_{i+1})$ w.r.t. $\idR=\idR(\mu)$ and that of $S(G_{i},G_{i+1})$ w.r.t. $\idR(\mu')$ give the same list of coefficients.  Thus, by the inductive hypothesis, we know that   all the parameters  that appear in $F'_i$, $i<r$, are eliminable, except (at most) $(r+1)r-2t$  of them.

 Finally  we can see that all the coefficients appearing in $F_r$ is eliminable, except at most   the $(r+1)$ coefficients $c_{rjk}$ such that  $x^{j+1}y^k  \in\idR$.
 
 Thus we can conclude that $\ed(\idR)\leq (r+1)+(r+1)r-2t+(r+1)=2\mu$.

\textbf{Step 2}: the case $n\geq 3$. To prove this general case, we first observe   that, as \texttt{RevLex} is a graded term ordering,   in the case $n=2$ the homogeneous stratum $\St(\idR(\mu))$ is canonically isomorphic to the  (non-homogeneous) stratum $\St(\idR_{af}(\mu))$ of the \lq\lq affine\rq\rq\ ideal in $k[x,y]$ having the same monomial basis as $\idR(\mu)$.  Thus, in order to construct $L(\idR(\mu,n))$ we can consider polynomials $F_i=x^{r-i}y^{i+\delta_i}+\sum C_{ijk}x^jy^k$ where $C_{ijk}=x^jy^j\sum_\gamma c_{ijk\gamma}Z^\gamma$, $Z^\gamma$ being any monomial in the variables $z_1, \dots, z_{n-1}$ of degree $r+\delta_i-j-k$.  Applying Procedure \ref{procedura}, we can consider each $C_{ijk}$ as   a variable. Thanks to what proved in Step 1 we know that   all the variables $C_{ijk}$ are eliminable except $2\mu$ of them and that there are no  relations among the remaining  ones. Finally we come back to the variables $c_{ijk\gamma}$: the non-eliminable among them are those appearing in the non-eliminable $C_{ijk}$.  We can compute their number following the same line of the proof of Step 1, where the non eliminable variables are explicitly listed. We observe for instance  that the formula of the dimension is quadratic in the number $n$ because for every non-eliminable $C_{ijk}$ the difference between the degree of $x^{r-i}y^{i+\delta_i}$ and that of $x^jy^k$, that is the degree of every $Z^\gamma$, is at most $2$.
\end{proof}

\bibliographystyle{plain}

\bigskip

\noindent Margherita Roggero, Lea Terracini \\
Dipartimento di Matematica dell'Universit\`{a} \\
          Via Carlo Alberto 10, 
           10123 Torino, Italy \\ 
 {\small margherita.roggero@unito.it } \\
{\small lea.terracini@unito.it } 
\end{document}